\documentclass[12pt]{amsart}
\usepackage{amssymb,amsthm}
\usepackage{eucal,mathrsfs,amsfonts}
\usepackage{color}
\usepackage{ulem,marginnote}
\usepackage[all]{xy}
\usepackage{fullpage}

\newtheoremstyle{mythm}%
  {\topskip}
  {\topskip}
  {\itshape}
  {}
  {\bfseries}
  {\\}
  {\parindent}
  {\thmname{#1}\thmnumber{ #2}\thmnote{ (\textit{#3})}}
\newtheoremstyle{mydef}%
  {\topskip}
  {\topskip}
  {}
  {}
  {\bfseries}
  {\\}
  {\parindent}
  {\thmname{#1}\thmnumber{ #2}\thmnote{ #3}}

\theoremstyle{mythm}
\newtheorem{lem}{Lemma}[section]
\newtheorem{thm}[lem]{Theorem}
\newtheorem*{thm*}{Theorem}

\newtheorem{cor}[lem]{Corollary}
\newtheorem{prop}[lem]{Proposition}

\theoremstyle{mydef}
\newtheorem{defn}[lem]{Definition}
\newtheorem{rem}[lem]{Remark}

\numberwithin{equation}{section}

\newcommand{\mbb}[1]{\mathbb #1}

\newcommand{\mc}[1]{\mathcal #1}

\newcommand{\oper}[1]{\operatorname{#1}}

\newcommand{\Gm}{\mbb G_m}

\newcommand{\Br}{\oper{Br}}

\renewcommand{\deg}{\oper{deg}}
\newcommand{\ram}{\oper{ram}}

\newcommand{\cha}{\oper{char}}

\newcommand{\Spec}{\oper{Spec}}

\newcommand{\calO}{\ensuremath{\mathcal{O}}}
\newcommand{\m}{\ensuremath{\mathfrak{m}}}
\newcommand{\slashfrac}[2]{\ensuremath{\raise1ex\hbox{#1}\kern-.2em/\kern-.30em\lower1ex\hbox{#2}}}

\newcommand{\res}{\ensuremath{\textrm{res}}}

\newcommand{\Q}{\mbb Q}
\newcommand{\Z}{\mbb Z}

\newcommand{\PP}{\mbb P}

\newcommand{\Jac}{\operatorname{Jac}}
\newcommand{\Pic}{\operatorname{Pic}}

\newcommand{\ov}{\overline}

\usepackage[OT2,T1]{fontenc}
\DeclareSymbolFont{cyrletters}{OT2}{wncyr}{m}{n}
\DeclareMathSymbol{\Sha}{\mathalpha}{cyrletters}{"58}

\def\<{\left<}
\def\>{\right>}

\title[Distinguishing algebras]{Distinguishing division algebras by finite splitting fields}

\author{Daniel Krashen}
\author{Kelly McKinnie}

\date{}

\begin{document}

\begin{abstract}
This paper is concerned with the problem of determining the number of division algebras which share the same
collection of finite splitting fields.  As a corollary we are able to determine when two central division
algebras may be distinguished by their finite splitting fields over certain fields.
\end{abstract}

\maketitle

\section{Introduction}

A major theme in the study of finite dimensional division algebras is
determining those field extensions of the center which split the algebra
--- i.e. such that the algebra becomes isomorphic to a matrix algebra when
the scalars are extended to this field extension. Despite the fact that
this is one of the major tools used to determine structural information
about such algebras, there are still a large number of open questions. In
this paper, we examine how much information is given by the finite
splitting fields of a central division algebra. To answer this, we
determine in certain cases how many distinct division algebras can share
the same collection of finite splitting fields, for example, showing that
any pair of quaternion algebras over the field $\mbb Q(t)$ which share the
same splitting fields must in fact be isomorphic.  This particular fact
answers a question originally posed to us by Peter Clark, and was our
original motivation this line of inquiry. Our paper further generalizes this to show, for example, that given
a division algebra $D$ over $\mbb Q(t)$ of prime period $p$, the
collection of division algebras of period $p$ sharing the same splitting field is always finite.

Independent parallel work of Garibaldi-Saltman \cite{SalGar} and Rapinchuk-Rapinchuk \cite{RapRap} has also
given the above result on quaternion algebras over the field $\mbb Q(t)$. Besides the fact that the methods we
use are quite distinct from these other two approaches, the results in these papers differ from ours in two
basic ways: First, \cite{SalGar} is concerned only with quaternion algebras (and symbols in higher cohomology
groups), and \cite{RapRap} is concerned only with period $2$ division algebras. Second, they are interested in
only maximal subfields, as opposed to finite splitting fields. By focusing on all finite splitting fields, as
opposed to simply the maximal subfields, we are able to prove results for more general fields than those
arising in \cite{RapRap}, as well as make statements concerning splitting fields for algebras of periods other
than $2$.

In the course of this paper, we also introduce some new techniques for
working with unramified Brauer classes on curves with rational points, 
for example, showing that such classes always arise as pullbacks of Brauer
classes on the Jacobian of the curve via the Abel-Jacobi map
(Proposition~\ref{surjective}).

\section{Statement of main results}

\begin{defn}
Let $k$ be a field and let $\alpha, \beta \in \Br(k)$. We write $\alpha
\equiv \beta$ if for every finite field extension $\ell/k$ we have
$\alpha_\ell$ is split if and only if $\beta_\ell$ is split. This defines
an equivalence relation on the elements of $\Br(k)$ and we let
$\|\alpha\|$ denote the equivalence class of $\alpha$.
\end{defn}

If $\|\alpha\| \ne \|\beta \|$ for some $\alpha,\,\beta \in \Br(k)$, then
there exists a finite field extension $\ell/k$ such that $\ell$ splits one
of $\alpha$ or $\beta$, but not the other.  Therefore, if $\|\alpha\| \ne
\|\beta\|$ we say that $\alpha$ and $\beta$ can be distinguished via
splitting fields.  Our main result is to say that if one knows a bound for
the order of the sets $\| \alpha\|$ over a field $k$, one may also
understand the potential size of these sets over $k(t)$.

\begin{thm} \label{main theorem 1}
Let $k$ be a field and $p$ a prime integer not equal to the characteristic
of $k$.  Let $\alpha \in \Br(k(t))[p]$ and suppose that the class $\alpha$
is ramified at $r$ distinct closed points. Then
\begin{enumerate}
\item if all $\alpha_0 \in
\Br(k)[p]$ satisfy $\#\|\alpha_0\| \leq N$ for some integer $N$, then
\[\|\alpha\| \leq N (p-1)^r. \]
\item if all $\alpha_0 \in \Br(k)[p]$ satisfy $\#\|\alpha_0\| < \infty$, then $\#\|\alpha\| < \infty$.
\end{enumerate}
\end{thm}

\begin{rem} \label{remark global etc} If $k$ is a higher local field, for example, an iterated Laurent series field
$k_0((t_1))\cdots((t_m))$ where $k_0$ is finite, local, or algebraically closed, then one may show that
the $p$-torsion part of the Brauer group, $\Br(k)[p]$, is finite for any $p$, and in particular, the
hypotheses of Theorem~\ref{main theorem 1}(1) will automatically hold for $N = |\Br(k)[p]|$. A description of
when $\Br(k)[2]$ is finite is given in \cite{Efrat:FiniteBrauer} in the case the characteristic of $k$ is not
$2$. We note also that the weaker conditions of Theorem~\ref{main theorem 1}(2) hold in the case that $k$
is a global field.
\end{rem}

\begin{cor} \label{split period 2}
Let $k$ be a field of characteristic not $2$, such that for all $\alpha_0 \in \Br(k)[2]$, we have
$\#\|\alpha_0\| =1$. Then for all $\alpha \in \Br(k(t))[2]$, $\#\|\alpha\| = 1$.
\end{cor}

In particular, this show that all the $2$-torsion elements in $\Br(\Q(t))$ may be distinguished via their
finite splitting fields, a question originally posed to us by Peter Clark.  This theorem is a consequence of
the proceeding results, for which we use the following notation:

For a closed point $x$ in a curve $X$, $\kappa(x)$ denotes the residue field, and for a Brauer class $\alpha
\in k(X)$ unramified at $x$, we let $\alpha|_x$ denote the specialization of $\alpha$ to the closed point $x$.
One may intepret this concretely by choosing $A$ to be a Azumaya algebra over the local ring $\mc O_{X,x}$
such that $A \otimes_{\mc O_{X,x}} k(X)$ represents the class $\alpha$, and then defining $\alpha|_x$ to be
the class of $A \otimes_{\mc O_{X,x}} \kappa(x)$.

The proof of Theorem \ref {main theorem 1} will depend on the following two theorems.
\begin{thm} \label{ramification theorem}
Let $k$ be a field and $p$ a prime integer not equal to the characteristic of $k$. Suppose that $\alpha, \beta
\in \Br(k(t))$ such that $\alpha \equiv \beta$. Then for every closed
point $x \in \mbb P^1(k)$, if we write $\ram_x(\alpha) = (L/\kappa(x), \sigma)$ and
$\ram_x(\beta) = (M/\kappa(x), \tau)$ then $L \cong M$ as extensions of
$\kappa(x)$.
\end{thm}
\begin{thm}[Distinguish classes with distinguishable specializations] \label{easycase}
Suppose $\alpha,\, \beta \in \Br(k(t))$ with $\alpha \neq \beta$, and
suppose there is a rational point $x \in (\PP^1_k)^{(1)}$ such that
$\ram_x\alpha=\ram_x\beta=0$ and $\beta|_x \not\equiv \alpha|_x$.  Then
$\beta \not\equiv \alpha$.
\end{thm}

The main content of the remainder of the sections of the paper will be to prove Theorem~\ref{ramification
theorem} (on page~\pageref{ramification theorem proof}) and Theorem~\ref{easycase} (on page~\pageref{easycase
proof}). Before doing so, we first illustrate how these theorems may be used to prove Theorem~\ref{main
theorem 1}.

\begin{proof}[Proof of Theorem \ref{main theorem 1}] Let $\alpha,\beta \in \Br(k(t))[p]$.  For any closed
point $x \in \mbb P^1_k$, let $\ram_x(\alpha) = (L_x/\kappa(x), \sigma)$.  By Theorem \ref{ramification
theorem}, if $\beta \in \|\alpha\|$, then for every closed point $x \in X$, if $\ram_x(\beta) = (M/\kappa(x),
\tau)$ then $M \cong L_x$ as extensions of $\kappa(x)$.  Therefore, $\tau=\sigma^i$ for some $1 \leq i \leq
p-1$ and $\alpha$ and $\beta$ ramifiy at the same set of closed points.  Let $\{x_j\}$ be the set of closed
points at which $\alpha$ ramifies.  For each of the $(p-1)^r$ possible sequences $(i_1,\ldots, i_r)$ with
$1\leq i_j\leq p-1$, let 
\[\|\alpha\|_{(i_1,\ldots, i_r)}=\{\beta \in \|\alpha\| \,|\,\ram_{x_j}\beta=(L_x/k(x),\sigma^{i_j}) \textrm{
for all }1\leq j \leq r\}\]
Then, $\|\alpha\|=\bigcup\|\alpha\|_{(i_1,\ldots,i_r)}$ with the union taken over all possible $(p-1)^r$
sequences. To prove part (1), it is only left to show that $\#\|\alpha\|_{(i_1,\ldots,i_r)}\leq N$ and for
this we use Theorem \ref{easycase}. 

In the case that $k$ is infinite, we may choose $x\in \PP^1(k)$ such that
$\ram_x\alpha=\ram_x\beta=0$. We show that $\#\|\alpha\|_{(i_1,\ldots,i_r)} \leq \#\|\alpha|_x\|$.  Note first
that $|_x : \|\alpha\|_{(i_1, \ldots, i_r)} \to \Br(k)$ is injective.  This follows since any two elements
$\beta_1,\beta_2 \in \|\alpha\|_{(i_1,\ldots,i_r)}$ have the exact same ramification sequence and therefore,
by the Auslander-Brummer-Faddeev sequence, $\beta_1 = \beta_2+\gamma$ for a constant class $\gamma \in
\Br(k)$. If $\beta_1|_x=\beta_2|_x$ then $\gamma|_x=\gamma $ is trivial, implying that $\beta_1=\beta_2$.  
In the case that $k$ is finite, it follows immediately from the Auslander-Brummer-Faddeev sequence (see e.g.,
\cite[6.9.3]{GiSz}) that
$\#\|\alpha\|_{(i_1,\ldots,i_r)} = 1$. 

Assume by way of contradiction that $\#\|\alpha\|_{(i_1,\ldots,i_r)} > \#\|\alpha|_x\|$.  Then, since the
specialization map $|_x$ is injective, $\beta|_x \notin \|\alpha|_x\|$ for some $\beta \in
\|\alpha\|_{(i_1,\ldots, i_r)}$.  By Theorem \ref{easycase} we can distinguish between $\alpha$ and $\beta$
using finite dimensional splitting fields, that is, $\beta \notin \|\alpha\|_{(i_1,\ldots, i_r)}$, a
contradiction.  Therefore, $\#\|\alpha\|_{(i_1,\ldots,i_r)}\leq N$ and $\#\|\alpha\|\leq (p-1)^rN$.  

To prove part (2) we use the terminology from above and set 
\[M=\max\{\#\|\alpha\|_{(i_1,\ldots,i_r)}\}.\]  
where the maximum is taken over all possible $(p-1)^r$ sequences.  As stated above, for any rational point $x
\notin \{x_j\}$, $\#\|\alpha\|_{(i_1,\ldots,i_r)}\leq\#\|\alpha|_x\|<\infty$, so $M$ is a finite number.  Then
$\#\|\alpha\|=\sum \|\alpha\|_{(i_1,\ldots,i_r)} \leq (p-1)^rM<\infty$.
\end{proof}

\section{Distinguishing Brauer classes via branched covers}

\begin{lem} \label{lem1}
Let $\phi:Y \to X$ be a branched cover such that there exists a
$k$-rational point $y \in Y(k)$ and let $x = \phi(y)$.  Then,
\begin{enumerate}
\item For any non-trivial constant class $\beta \in \Br(k)\subseteq
\Br(k(X))$, $\beta_{k(Y)} \ne 0$.
\item If $\phi:Y \to X$ is unramified at $y$ then for any class $\alpha
\in \Br(k(X))$, $\ram_y\alpha_{k(Y)}=\ram_x\alpha$.
\end{enumerate}
\end{lem}
\begin{proof} For part (1), consider the commutative diagram
\begin{equation} \label{d1}
\xymatrix{
\Br(k) \ar@{^{(}->}[r] \ar[d]_{\textrm{id}} & \Br(X) \ar[d]_{\phi^*}
\ar@{^{(}->}[r] & \Br(k(X)) \ar[d]^{\res} \ar[r]^(.45){\ram_x}
& H^1(k,\Q/\Z)\ar[d]^{e \cdot \res} \\
\Br(k) \ar@{^{(}->}[r] & \Br(Y) \ar@{^{(}->}[r] & \Br(k(Y))
\ar[r]^(.45){\ram_y} & H^1(k,\Q/\Z),
}
\end{equation}
where $e$ is the ramification index of $\phi$ at $y$.
The arrows $\Br(k) \hookrightarrow \Br(X)$ and $\Br(k)
\hookrightarrow\Br(Y)$ are injections because of the existence of
sections given by the rational points. Therefore, $k(Y)$ cannot split
$\beta$. Part (2) follows from the right hand side of the  commutative
diagram (\ref{d1}),(see \cite{Sal:LN}) and the fact that $e=1$ since $\phi:Y
\to X$ is unramified at $y$.\end{proof}

Although our goal is to prove statements about Brauer classes over $k(t)$, the function field of $\mbb P^1_k$,
other curves and their Jacobians will naturally arise in the process. We will therefore shift focus for a
while and consider the ramification and splitting behavior of Brauer classes of curves over more general
curves. 

Let $X$ be a smooth projective curve over $k$ containing a
rational point $x \in X(k)$.  Let $\phi:X \to \Jac X$ be the Albanese
map taking $x$ to $[0]$, the identity element.  Throughout this section
we will set $J=\Jac X$.  The next two lemmas collect facts about $X$
and $J$ which we use in Lemma \ref{l3} to produce a cover of $X$ which
makes Brauer classes on $X$ constant.

\subsection{Make unramified}
\begin{lem}[Make it unramified]\label{make it unramified}
Let $k$ be a field and $X$ a smooth projective $k$ curve.  Let $\alpha \in \Br(k(X))[m]$ with
$(m,\textrm{char}(k))=1$ and let $D$ be the ramification locus of $\alpha$.  Assume there exists $x\in X(k)$,
with $\ram_x(\alpha)=0$.  Then, there exists a branched cover $\psi:Y \to X$ and a rational point $y \in Y(k)$
such that
\begin{enumerate}
\item $\psi(y)=x$, and
\item $\alpha_{k(Y)}$ is unramified at all $y\in Y^{(1)}$ so that $\alpha_{k(Y)} \in \Br(Y)[m]$.
\end{enumerate}
\end{lem}

\begin{proof} 
Choose a closed point $b \neq
x$ in $X$, with $b$ not in the support of $D$. We wish to find a rational function on $X$ which is regular
away from the point $b$, which vanishes on the ramification locus $D$ with only simple zeroes, and which is
nonzero at the rational point $x$. Such a function would exactly correspond to a global section of $\mathscr
L(Nb - D)$ which is also not a section of $\mathscr L(Nb - D - x)$ or $\mathscr L(Nb - D - d)$ for any $d$ in
the support of $D$.

Choose $N > \left(2g - 2 + 2\deg(D)\right)$ where $g$ is the genus of $X$, which in particular ensures that 
\[\deg(Nb - D), \deg(Nb - D - x), \deg(Nb - D - d) > 2g - 2,\]
Therefore, $l(K-Nb-D-a)=l(K-Nb)=0$ where $K$ is the the class of the canonical divisor (see
\cite[IV.1.3.3]{Hartshorne}).  By Riemann-Roch, 
\begin{eqnarray*}
l(Nb-D-z)&=&\deg(Nb-D-z)+1-g\\
&<&\deg(Nb-D)+1-g\\
&=&l(Nb-D).
\end{eqnarray*}  
For $z = x$ or $z$ a closed point in the support of $D$. In particular, by choosing $N$ sufficiently large we
find that
\[\emptyset \neq \Gamma\left(\mathscr L(Nb - D)\right) \setminus \left( \Gamma \left(\mathscr L(Nb - D -
x)\right) \cup \bigcup_{d \in supp(D)} \Gamma
\left(\mathscr L(Nb - D - d) \right) \right) \]
That is, there is a global section $f$ of $\mathscr L(Nb-D)$ with the germ $f_x \notin \m_{X,x}$ and with
only simple zeros on $D$.  

Let $f(x)$ be the image of $f_x$ in $\calO_{X,x}/\m_{X,x}\cong k$.  Set 
\[L=k(X)[T]/(T^m-ff(x)^{m-1}).\]  
Since every closed point $d \in \textrm{supp}(D)$ has multiplicity 1, our choice of $f$ ensures $f \in \m_{X,d}-\m_{X,d}^2$ and further, since $f(x)\in k\subset \calO_{X,d}$ is a unit, $ff(x)^{m-1} \in \m_{X,d}-\m_{X,d}^2$.  Therefore, by Eisenstein's criterion, the polynomial $T^m-ff(x)^{m-1}\in \calO_{X,d}[T]$ is irreducible.  In particular, $L$ is a field and $T$ is integral over $\calO_{X,d}$ for all $d \in \textrm{supp}(D)$.

Let $Y$ be the normalization of $X$ in $L$ with morphism $\psi:Y \to X$.  We show that $Y$ satisfies the
condition of the lemma.  Let $p\in X$ be a closed point.  If $\ram_p\alpha=0$, then for any point $q\in Y$ lying
over $p$, $\ram_q\alpha_{k(Y)}=0$.  If $\ram_p\alpha \ne 0$ then $p\in\mathrm{supp}(D)$.  Let $q \in Y$
with $\psi(q)=p$.  Let $\pi_p$ be a uniformizer for $\calO_{X,p}$ so that $ff(x)^{m-1}=\pi_pu$ with $u$ a unit
in $\calO_{X,p}$.  Let $v_q$ be the valuation of $\calO_{Y,q}$ extending that of $\calO_{X,p}$.  Then
$v_q(T^m)=mv_q(T)=v_q(\pi_p)=1$.  Since $[L:k(X)]=m$, it follows that the point $q$ has
ramification index $e_q=m$.  Therefore, by the standard commutative diagram
\[\xymatrix{
\Br(k(X))\ar[r]^(.4){\ram_p}\ar[d]&H^1(k(p),\Z/m\Z)\ar[d]^{e_q}\\
\Br(L)\ar[r]^(.34){\ram_q}&H^1(k(p),\Z/m\Z)
}\]
$\ram_q \alpha_{L}=m\cdot\ram_p\alpha=0$.

It is only left to show that there exists a $k$-rational point in $Y$ lying over $x\in X(k)$.  To do this we
compute the fiber $Y_x$ of $Y$ over $x$.  First we note that for any $p\notin \mathrm{supp}(D)$, $ff(x)^{m-1}
\in \calO_{X,p}-\m_{X,p}$ is a unit.  Moreover, by assumption $m$ is a unit in $k\subset \calO_{X,p}$.
Therefore the morphism $\Spec\left(\calO_{X,p}[T]/(T^m-ff(x)^{m-1})\right) \to \Spec \calO_{X,p}$ is \'etale
and in particular, $\calO_{X,p}[T]/(T^m-ff(x)^{m-1})$ is integrally closed and hence the integral closure
of $\calO_{X,p}$ in $L$. In particular,
\begin{equation*}
\begin{split}
Y_x=&\Spec\left((\textrm{Int. Clos.}_L\calO_{X,x})\otimes_{\calO_{X,x}}\frac{\calO_{X,x}}{\m_x}\right)\\
=&\Spec\left(\frac{k[T]}{T^m-f(x)^m}\right)\\
=&\Spec k \times \Spec\left(\frac{k[T]}{g(T)}\right)
\end{split}
\end{equation*}
where $T^m-f(x)^m=(T-f(x))g(T)$ and $g(f(x))\ne 0$.  This shows that $Y_x$ has a $k$-rational point and hence that there is a $k$-rational point lying over $x$.\end{proof}

\subsection{Make constant}

\begin{lem}\label{Picard-Brauer sequence}
Let $V$ be a smooth projective $k$-variety containing a rational point
$v \in V(k)$. Then we have a short exact sequence 
\[0 \to \Br(k) \to \Br(\ov{V}/V) \to H^1(k, \Pic \ov V) \to 0 \]
arising from the $E_2$-terms of the Hochschild-Serre spectral sequence
for $V$.
\end{lem}
\begin{proof}
See for example, \cite{Skor:TRP}~Corollary 2.3.9. If we write
\[E^2_{p,q} = H^p(k, H^q(\ov{V}, \mbb G_m) \implies H^{p+q}(V, \mbb G_m) = E_{p+q} \]
then the above sequence may be identified with the terms in the spectral sequence as:
\[0 \to E_{2, 0}^2 \to \left[ \ker \left(E_2 \to E_{0,2}^2\right)\right] \to E_{1,1}^2 \to 0 \]
\end{proof}

\begin{prop}[$\Br(J)$ surjects onto $\Br(X)$] 
The pullback map $\phi^*:\Br(\ov J/J) \to \Br(X)$ is surjective.
\label{surjective}
\end{prop}
\begin{proof} 
Using the fact that pullback of cohomology classes induces a morphism
of Hochschild-Serre spectral sequences
\[\xymatrix{
H^p(k, H^q(\ov J, \Gm)) \ar[d] \ar@{=>}[r] & H^{p+q}(J, \Gm) \ar[d] \\
H^p(k, H^q(\ov V, \Gm)) \ar@{=>}[r] & H^{p+q}(V, \Gm) \\
}\]
combined with Lemma~\ref{Picard-Brauer sequence}, we see that this
morphism of spectral sequences yields a morphism of short exact
sequences:
\begin{equation}\label{picard-brauer diagram}
\xymatrix{
0 \ar[r] & \Br(k) \ar[d]^{id} \ar[r] & \Br(\ov{J}/J) \ar[d] \ar[r] &
H^1(k, \Pic \ov J) \ar[r] \ar[d] & 0  \\
0 \ar[r] & \Br(k) \ar[r] & \Br(\ov X/X) \ar[r] & H^1(k, \Pic \ov X)
\ar[r] & 0 
}
\end{equation}

From the short exact sequence 
\[0 \to \Pic^0X\to \Pic X \to \Z\to 0,\]
we obtain an isomorphism $H^1(k,\Pic^0\overline X) \cong H^1(k,\Pic
\overline X)$.  Since $\phi^*$ induces an isomorphism $\Pic_0 \overline
J \cong \Pic_0 \ov X$, we may use the commutative
diagram
\[\xymatrix{
H^1(k, \Pic^0 \overline J) \ar[d] \ar[r]^{\phi^*}_\cong & H^1(k, \Pic^0
\overline X) \ar[d]^\cong \\
H^1(k, \Pic \overline J)\ar[r]^{\phi^*}&H^1(k,\Pic \overline X)
}\]
to verify that the pullback map induces a surjection $H^1(k, \Pic \ov J)
\to H^1(k, \Pic_0 \ov X)$.

By Tsen's Theorem, $\Br(k(\overline X))=0$. Therefore, $\Br(\overline
X)\subseteq \Br(k(\overline X))$ is also trivial (\cite{Milne:EC},
IV2.6), and in particular, we may identify $\Br(\ov{X}/X) = \Br(X)$.

It now follows from a diagram chase using diagram~(\ref{picard-brauer
diagram}) that the pullback map $\Br(\ov J/J) \to \Br(X)$ is surjective.
\end{proof}

For any integer $n$, let $\widetilde X_n$ be the fiber product
\begin{equation}\xymatrix{
\widetilde X_n \ar[r]^q\ar[d]_p&J\ar[d]^n\\
X\ar[r]_{\phi}&J
}\label{d2}\end{equation}
where $n$ is the multiplication by $n$ map and as always, $J=\Jac X$ and
$\phi:X \to J$ is the Albanese map given by the rational point $x \in
X(k)$. Since the multiplication by $n$ map is \'etale, $p$ is an \'etale covering.

\begin{lem}[Pulling back to $\widetilde X_n$.]  Let $\widetilde X_n$ be given as above.
\begin{enumerate}
\item Let $\alpha \in \Br(X)[m]$.  Then $p^*\alpha\in \Br(\widetilde X_{2m})$ is a constant class.
\item Let $x\in X(k)$ be a rational point.  Then for any $n$, there exists a $k$-rational point $\tilde x \in \widetilde X_n(k)$ with $p(\tilde x)=x$.
\end{enumerate}
\label{l3}
\end{lem}
\begin{proof} (1) By the surjectivity of $\phi^*$ given in Proposition~\ref{surjective} there exists $\beta \in \Br(J)$ such that $\phi^*(\beta)=\alpha$.  Let $\bar\beta\in H^1(k,\Pic \overline J)$ be the image of $\beta$.  To show that $(2m)^*(\beta)\in \Br(J)$ is constant it is enough to show that $\overline{(2m)^*\beta}= 0$ in $H^1(k,\Pic \overline J)$.  This will follow from the commutativity of
\[\xymatrix{
\frac{\Br(\overline{J}/J)}{\Br(k)}\ar[r]^(.4){\cong}\ar[d]_{n^*}&H^1(k,\Pic \overline J)\ar[d]^{n^*}\\
\frac{\Br(\overline{J}/J)}{\Br(k)}\ar[r]^(.4){\cong}&H^1(k,\Pic\,\overline J)
}\]
for all $n \in \Z$, which follows from the naturality of the Hochschild-Serre spectral sequence.  Here the map $n^*$ on the right hand side is the map on cohomology gotten from $n^*:\Pic\,\overline J \to \Pic\,\overline J$, the map pulling back line bundles. 

By e.g., 
\cite[pg 59, Corollary 3]{Mum:AV}, for a line bundle $\mathcal L$ on $J$, $n^*(\mathcal L)=\mathcal L^{\frac{n^2+n}{2}}\otimes (-1)^*\mathcal L^{\frac{n^2-n}{2}}$.  Notice that the map 
\begin{equation*}
\begin{split}
(2m)^*:\Pic\, J &\to \Pic\,J,\\ 
\mathcal L &\mapsto \mathcal L^{m(2m+1)}\otimes (-1)^*\mathcal L^{m(2m-1)}
\end{split}\end{equation*}
factors as $(2m)^*=f\circ g$ where $g(\mathcal L) = \mathcal L^m$ and $f(\mathcal L) = \mathcal L^{2m+1}
\otimes (-1)^*(\mathcal L^{2m-1})$.  Consequently $(2m)^*:H^1(k,\Pic\,\overline J) \to H^1(k,\Pic \overline
J)$ also factors as $f\circ g$, where $g$ is multiplication by $m$.  Therefore, since $\bar\beta$ has order
dividing $m$ (since $\alpha \in \Br(X)[m]$), $(2m)^*(\bar\beta)=0$.  In other words, $(2m)^*(\beta)$ is a
constant class.  Using the notation from diagram (\ref{d2}), $p^*\alpha=q^*(2m)^*\beta$, thus we see
$p^*\alpha$ is also a constant class.

(2)  By the definition of $\phi$, $\phi(x)=[0]\in J$ where $x$ is our rational point.  Therefore, by the
definition of the fiber product, to produce a rational point of $\widetilde X_n(k)$, we need only produce a
rational point of $J$ which maps to $[0]$ under $n$.  We simply take the origin itself, since $n[0]=[0]$ for
all $n \in \Z$.
\end{proof}

For simplicity we state an immediate corollary.
\begin{cor}[Make unramified classes constant]  Let $\alpha \in \Br(X)[m]$.  Then there is an \'etale cover of
curves $p:\widetilde X \to X$ and a rational point $\tilde x \in \widetilde X(k)$ so that $p(\tilde x)=x$ and
$p^*(\alpha) \in \Br(k)$ is a constant class.
\label{from unramified to constant}
\end{cor}

\begin{lem}[Make it constant] \label{make it constant}
Let $X/k$ be a smooth projective curve with $x \in X(k)$, and $\alpha \in \Br(k(X))[m]$ with $\gcd(m,\textrm{char}\,k)=1$ and $\ram_x(\alpha) = 0$.  There exists a branched cover $\phi: Y \to X$ such that 
\begin{enumerate}
\item $\alpha_{k(Y)}$ is constant,
\item there is a $y \in Y(k)$ such that $\phi$ is unramified at $y$ and $x = \phi(y)$.
\end{enumerate}
\end{lem}
\begin{proof}
By  Lemma \ref{make it unramified} there exists a branched cover $\psi:Y \to X$ and a rational point $y \in
Y(k)$ so that $\psi(y)=x$ and $\alpha_{k(Y)} \in \Br(Y)[m]$.  So, without loss of generality, we may assume
$\alpha \in \Br(X)[m]$, where $X$ is a smooth projective curve over $k$ with rational point $x$.  We are
then done by Corollary \ref{from unramified to constant}
\end{proof}

We may now complete the proof of Theorem~\ref{easycase}:
\begin{thm*}[\ref{easycase}]
Suppose $\alpha,\, \beta \in \Br(k(t))$ with $\alpha \neq \beta$, and
suppose there is a rational point $x \in (\PP^1_k)^{(1)}$ such that
$\ram_x\alpha=\ram_x\beta=0$ and $\beta|_x \not\equiv \alpha|_x$.  Then
$\beta \not\equiv \alpha$.
\end{thm*}

\begin{proof} \label{easycase proof} Using Lemma~\ref{make it constant}, once for the class $\alpha$ and then
again for the pullback of the class $\beta$ to the resulting cover, we may find a branched cover $\phi : Y \to
\PP^1_k$ such that both $\alpha_{k(Y)}$ and $\beta_{k(Y)}$ are constant classes and such that there is a
$k$-rational point $y \in Y(k)$ with $\phi(y) = x$.  Then $\alpha_{k(Y)}|_y=\alpha|_x$ and
$\beta_{k(Y)}|_y=\beta|_x$ and therefore, $\alpha_{k(Y)}=(\alpha|_x)_{k(Y)}$ and
$\beta_{k(Y)}=(\beta|_x)_{k(Y)}$.

By hypothesis $\beta|_x \not\equiv \alpha|_x$ and therefore, modulo
switching $\alpha$ and $\beta$, there exists a finite field extension
$L/k$ such that $(\alpha|_x)_L =0\ne(\beta|_x)_L$.  Since $Y$ has a
rational point, $\Br(L) \to \Br(L(Y))$ is injective and so it follows
that we have $0=\alpha_{L(Y)} \ne \beta_{L(Y)}$.\end{proof}

\section{Distinguishing ramified classes from constant classes}

\begin{lem}[Make ramification occur at a rational point]\label{make
ramification rational}
Let $\alpha\in\Br(k(t))[m]$ where $m$ is not divisible by $p = \cha{k}$
and assume $\ram_x\alpha=(L/k(x),\sigma)$ for some closed point $x \in 
X$.  Then there is a
$k(x)$-rational point $x'$ in $X_{k(x)}$ lying over $x$ and, for any
such point $x'$, we have $\ram_{x'}\alpha = (L/k(x),\sigma^{p^n})$ for
some $n$.
\end{lem}
\begin{proof} 
Without loss of generality, we may assume that $\ram_x \alpha \neq 0$.
Let $k' = k(x)$ be the residue field of the point $x$, and let $x' \in
\PP^1_{k'}$ be a $k'$ rational point lying over $x$. 
The lemma follows from the standard commutative diagram
\[\xymatrix{
\Br(k(t))' \ar[r]^{\ram_x}\ar[d]_{\mathrm{res}}&H^1(k',\Q/\Z)'\ar[d]^{e}\\
\Br(k'(t))'\ar[r]^{\ram_{x'}}&H^1(k',\Q/\Z)'
}\]
and the fact that $e=[k':k]_i$, is
the inseparability degree of $k'/k$  see e.g. \cite[9.19]{GMS}. 
\end{proof}

\begin{lem}[Distinguish between ramified and constant]\label{ram constant}
Let $p$ be a prime with $(\cha k,p)=1$ and let $\alpha,\beta \in \Br(k(X))[p]$ be classes such that
\begin{enumerate}
\item $\alpha$ is ramified at a rational point $x \in X(k)$ and,
\item $\beta$ is a constant class (i.e., $\beta \in \Br(k)[p]\subseteq\Br(k(X))$).
\end{enumerate}
Then $\alpha \not\equiv \beta$.
$\beta_L$ is split and the other is not.
\label{hardcase}
\end{lem}

\begin{proof} Since neither $\alpha$, $\beta$, nor the ramification of $\alpha$
at $x$ will be split by a prime to $p$ extension, we may assume without
loss of generality that $\mu_p\subset k^*$.  Let $\omega$ be a primitive
$p$-th root of unity.  Let $t\in k(X)$ be a uniformizer for the local
ring $\mathcal{O}_{X,x}$ at $x$.  Assume to begin with that
$\alpha=(a,t)$ is a symbol with $a\in k^*$ a representative of the
ramification $\ram_x\alpha \in H^1(k,\Z/p\Z)\cong k^*/(k^*)^p$.  Let $Y$
be a model of the extension $L = k(X)[s]/(s^p + ta^{-1})$ over $X$ (we can take, $Y$ to be the normalization
of $X$ in $L$).  Then one may check that since $\mathcal O_{X,x} [s]/(s^p + ta^{-1})$ is integrally closed,
$Y$ has a
$k$-rational point lying over $x$ (at $s = t = 0$). Futher, we note that $u = t/s$ satisfies $u^p = -at^{p-1}$,
and so the extension $L$ splits $\alpha$
since $(a,t)_{L}=(a, au^p)_L = (a,-at^p)_L=(a,-a)_L=0$.  By Lemma \ref{lem1}(1), since
$Y$ has a rational point, $L=k(Y)$ does not split $\beta$.

We now consider the general case $\alpha \in \Br(k(X))[p]$ with
ramification at $x$.  Write $\alpha=\gamma+(a,t)$ where $\gamma$ is
unramified at $x$ and as above $a \in k^*$ satisfies $\bar
a=\ram_x\alpha$.  By Lemma~\ref{make it constant} there exists a
branched cover $Y \to X$ with a rational point $y \in Y(k)$ lying over
$x$ so that \begin{enumerate}
\item $\gamma_{k(Y)}$ is constant,
\item $\ram_q\alpha_{k(Y)}=\ram_p\alpha=\bar a$, and
\item $\beta$ is not split by $k(Y)$ because $Y$ has a rational point, but is still (of course) a constant class.
\end{enumerate}
So we may assume without loss of generality that $\alpha=\gamma+(a,t)\in
\Br(k(X))[p]$ where $\gamma$ is a constant class. Note that since the
cover $Y \to X$ is unramified at $y$ it follows that $t$ pulls back to a
uniformizer at $y$. By the Merkurjev-Suslin Theorem, we may fix a representation for $\gamma$ as a sum of
symbols, $\gamma = \sum_{i=1}^n(a_i,b_i) \in H^2(k,\mu_p)$.

Suppose that for some $i$, the subgroups of $k^*/(k^*)^p$ generated by $a_i$ and $a$ do not coincide. We claim
that the ramification of the symbol $(a, t)$ is not split by the
extention $\ell_i = k(\sqrt[p]{a_i})$. Since the cover $X_{\ell_i} \to
X$ is unramified, we obtain a commutative diagram
\[\xymatrix{
\Br(k(X)) \ar[r]^{ram} \ar[d]_{res} & H^1(k, \mu_p) \ar[d]^{res} \\
\Br(\ell_i(X_{\ell_i})) \ar[r]_{ram} & H^1(\ell_i, \mu_p)
}\]
In particular, by Kummer theory, it follows that the extensions $\ell_i$
and $k(\sqrt[p]a)$ are linearly disjoint, implying that the ramification
of $(a, t)$ remains nontrivial at the $\ell_i$-rational point lying over
$x$. We proceed to extend scalars to each such $\ell_i$ successively
(and replacing $k$ by $\ell_i$).  This cannot split the ramification of
$(a, t)$ and hence will not split the class $\alpha$. If $\beta$ becomes
split, then we have shown $\alpha \not\equiv \beta$ and so we are done.

After extending scalars as above, we may assume that $\beta$ is not split and that for each $(a_i, b_i)$ in
the decomposition of $\gamma$, $a_i$ has the same $p$-power residue class as $a^r$ for some $r$.
In particular, we find $(a_i, b_i) = (a^r, b) = (a, b^r)$, and so
without loss of generality, we may combine all such symbols to write
\[ \alpha = (a, t) + (a, b)  = (a, bt) \]
for some $b \in k^*$.  
However in this case we have reduced to the original one considered
since $bt$ is another uniformizer for the rational point $x\in X(k)$.
That is, $\alpha \not\equiv \beta$ since there exists a finite extension
$L/k(X)$ which splits $\alpha$ without splitting $\beta$.\end{proof}

We may now complete the proof of Theorem~\ref{ramification theorem}

\begin{thm*}[\ref{ramification theorem}]
Let $k$ be a field and $p$ a prime integer. Suppose that $\alpha, \beta
\in \Br(k(t))$ such that $\alpha \equiv \beta$. Then for every closed
point $x \in X$, if we write $\ram_x(\alpha) = (L/\kappa(x), \sigma)$ and
$\ram_x(\beta) = (M/\kappa(x), \tau)$ then $L \cong M$ as extensions of
$\kappa(x)$.
\end{thm*}
\begin{proof} \label{ramification theorem proof}
By Lemma~\ref{make ramification rational}, we may replace $k$ by a finite extension, and assume that $x \in
\mbb P^1(k)$. Suppose that $M \not\cong L$ as extensions of $k$. Since they are degree $p$ extensions, it
follows that they are linearly disjoint and so by extending scalars from $k$ to $M$ we may assume that in fact
$\beta$ is unramified at $x$ and $\alpha$ is ramified at $x$. By Lemma~\ref{make it constant}, we may find a
branched cover $\phi : X \to \mbb P^1_k$ with a rational point $x'$ lying over $x$, and such that $\beta_{k(X)}$ is
a constant class and $\phi$ is unramified at $x'$. By Lemma~\ref{lem1}, it follows that $\alpha_{k(X)}$ is
ramified at $x'$. Finally, we may conclude from Lemma~\ref{ram constant} that $\alpha \not \equiv \beta$ as
desired.
\end{proof}

\bibliographystyle{alpha}
\bibliography{citations}

\end{document}